\documentclass[12pt,a4paper,twoside]{article}
\newtheorem{theorem}{Theorem}

\newtheorem{proposition}{Proposition}

\def\BbbR{{I \! \! R}}
\def\BbbP{{I \! \! P}}
\def\BbbE{{I \! \! E}}

\def \R{\BbbR}

\def \P{\BbbP}
\def \E{\BbbE}

\def\1{{1  \! \! \hspace{0.22mm} \mbox{l}}}

\begin{document}
\title{A functional non-central limit theorem for jump-diffusions with periodic
coefficients driven by stable L\'{e}vy-noise\\[3ex]}
\author{Brice Franke\\ 
 \footnotesize  Fakult\"at f\"ur Mathematik, Ruhr Universit\"at Bochum\\[-1ex]
   \footnotesize  Universit\"atsstr. 150, D-44780 Bochum, Germany\\[-1ex] 
      \footnotesize  e-mail: Brice.Franke@ruhr-uni-bochum.de}

\maketitle

\abstract{We prove a functional non-central limit theorem for
jump-diffusions with periodic coefficients driven by stable
L\'{e}vy-processes with stability index $ \alpha>1 $. 
The limit process turns out to be an $ \alpha $-stable 
L\'{e}vy process with an averaged jump-measure. Unlike in the situation
where the diffusion is driven by Brownian motion, there is no drift related 
enhancement of diffusivity.}\\

{\it Keywords:} non-central limit theorem, homogenization,
asymptotic analysis, periodic diffusion, stable L\'{e}vy-process\\

{\bf MSC:} 60F05, 60G52 \\

\section{Introduction}

There are only few situations, where partial differential equations can be 
solved explicitly; for example, those with constant coefficients. It is 
therefore important to find tools to reduce a large class of equations to 
those of constant coefficients. Equations with periodic coefficients form such 
a class and have been analyzed by homogenization techniques. Indeed, it can be 
shown that suitable scalings of solutions to equations with periodic 
coefficients converge in an appropriate sense toward solutions to equations 
with constant coefficients. In the book of Bensoussan, Lions and Papanicolaou 
homogenization theorems for second order parabolic problems were proved 
with arguments from stochastic analysis (see \cite{[BenLioPap]} p.345).
They proved a central limit theorem for an It\^{o}-diffusion related to 
the original parabolic problem to obtain homogenization-results.
The key to their approach is the well known relation between the parabolic 
equation 
$$  \partial_tu(x,t)
 =\sum_{i,j,k}^da_{ij}(x)\partial_{x_i}\partial_{x_j}u(x,t)
    +\sum_{i=1}^db_i(x)\partial_{x_i}u(x,t) $$ and 
the stochastic differential equation 
$$ dZ_t=b(Z_t)dt+\sigma(Z_t)dB_t, $$
with $ a_{ij}=1/2\sum_k\sigma_{ik}\sigma_{jk} $.
If $ b $ and $ \sigma $ are assumed to be periodic and sufficiently smooth, 
then there exists an invariant measure $ \pi $ for the process $ Z $.
It was proved in \cite{[BenLioPap]} and with a different method in 
\cite{[Bha]}, that the renormalized processes
$$ Z^{(n)}_t:=
\sqrt{n}^{-1}\left(Z_{nt}-Z_0-nt\int_{[0,1]^d}bd\pi\right) $$
converge in distribution toward a Wiener-process with covariance-matrix given 
by 
$$  \Sigma_{\alpha\beta}=\int_{[0,1]^d}\sum_{i,j}
 (\delta_{i\alpha}-\partial_{x_i}\psi_\alpha)
  a_{ij}(\delta_{j\beta}-\partial_{x_j}\psi_\beta)d\pi,$$
where $ \psi_\beta $ is a periodic solution of the following Poisson-problem 
\begin{eqnarray*} \sum_{i,j}a_{ij}(x)\partial_{x_i}\partial_{x_j}\psi_\beta(x)
    +\sum_{i=1}^db_i(x)\partial_{x_i}\psi_\beta(x)
 =b_\beta(x)-\int_{[0,1]^d}b_\beta d\pi 
\end{eqnarray*}
(see \cite{[BenLioPap]} p.394). The emergence of an additional diffusion coefficient relateted to the 
Poisson-problem in the limit-process is sometimes called enhancement of diffusivity in the 
literature. In the book of Bensoussan, Lions and Papanicolaou the question aroused, 
whether a theorem of this type holds, if the Brownian motion in the stochastic 
differential equation for $ Z $ is substituted by stable L\'{e}vy noise 
(see \cite{[BenLioPap]} p.531). We will prove that under suitable 
conditions a theorem of this kind is valid. However we will have to restrain ourselves to the case where
the stability index $ \alpha $ of the stable L\'{e}vy-noise is larger than one. 
This reflects the fact that the part of the generator related to the drift has differentiation order
equal to one whereas the part related to the noise has differentiation order equal to 
$ \alpha $. For $ \alpha\leq1 $ ellipticity, which is relevant for the existence of transition-densities and 
invariant probabilities, cannot hold. 
The difference from the case of Brownian motion is that the limiting stable L\'{e}vy 
process does not depend on the solution to an associated Poisson-problem. 
Indeed, it only depends on the drift through the invariant measure. This 
non-existence of diffusivity enhancement is due to the strong scaling $ n^{-1/\alpha} $,
which overscales the quadratic variation of an additional randomness comming from the drift.\\[1mm]
In a recent paper the author studied the behaviour of stable-like processes with periodic coefficients 
under scaling (see \cite{[Fra]}). Stable-like processes are generalisations of stable L\'{e}vy processes,
where the stability-index $ \alpha $ depends on the location of the diffusion in the state-space. Similar
to the situation in the present paper one can prove weak convergence toward a stable L\'{e}vy 
process under some conditions on the scale-function $ \alpha(x) $. However, since the jump-diffusions
considered in \cite{[Fra]} have no drift, the proof of the limit-theorem is less difficult.
In the present paper a considerable part of the effort in the main results proof goes into the control of the drift. Following an idea from \cite{[BenLioPap]} the It\^{o}-formula and the solution of an associated Poisson problem are used to
replace the drift in the SDE by a term which on large scales behaves like a martingale with bounded jumps. 
This martingale is the reason for the diffusivity enhancement in the Brownian motion case where $ \alpha=2 $. 
In our situation it will turn out that the scaling is to strong such that the additional martingale-term converges to zero 
in the scaling-procedure. This can be motivated by looking at the quadratic variation of the additional martingale.
The factor $ n^{-1/\alpha} $ in front of the martingale becomes a factor $ n^{-2/\alpha} $ on the level of  quadratic variation. However the growth of the quadratic variation of the martingale due to time contraction is
linear in $ n $. This explains the non-existence of drift-related diffusivity enhancement in the 
homogenization of SDE's with non-gaussian stable L\'{e}vy noise.

\section{Preliminaries}

Let $ (\Omega,{\cal F},\P,({\cal F}_t)_{t\geq0}) $ be a filtered
probability space satisfying the usual conditions, i.e. $ {\cal F}_t $ is 
right-continuous and complete (see \cite{[JacShi]} p.2). Further, let 
$ D_{\R^d}([0,\infty[) $ be the space of c\`{a}dl\`{a}g-functions with the 
usual Skorohod topology (see \cite{[JacShi]} p.289). For a finite 
Borel-measure $ \mu $ on $ S^{d-1} $ and 
$ \alpha\in]0,2[ $, we use spherical coordinates 
$ (\hat{y},|y|):=(y/|y|,|y|)\in S^{d-1}\times\R^+ $
to define a $ \sigma $-finite measure on $ \R^d\backslash\{0\} $ as follows
\begin{eqnarray*}
   \nu(dy):= \mu(d\hat{y})|y|^{-\alpha-1}d|y|.
\end{eqnarray*}
The measure $ \nu $ is the compensator of an $ \alpha $-stable 
$ {\cal F}_t $-adapted L\'{e}vy-process in $ \R^d $. 
We will denote by $ L $ its unique c\`{a}dl\`{a}g-modification.\\[2ex]
We proceed under the following assumptions:\\[2mm]
{\bf A1} We assume that $ 1<\alpha<2 $.\\[1mm]
{\bf A2} Let $ \mu $ be symmetric i.e.: $ \mu(A)=\mu(-A) $ for all
measurable $ A\subset S^{d-1} $. \\[1mm]
{\bf A3}  For suitable positive constants $ C_1, C_2 $ we have
$$ C_1\leq\int_{S^{d-1}}|\langle v,\varphi\rangle|^\alpha\mu(d\varphi)
  \leq C_2. $$ 
{\bf Remarks:} Some comments on our assumptions have to be made:\\[0.5mm]
1) The assumption A1 is necessary since for $ 0<\alpha\leq1 $ the first order 
operator associated to the drift would have an order of differentiation exceeding 
the order of differentiation of the generator of $ L $. 
Ellipticity could not hold anymore.\\[0.5mm]
2) Under the assumption A2 the center resp. drift of the L\'{e}vy-process 
$ L $ is equal to zero (see \cite{[Sat]} p.39 et sqq.).\\[0.5mm]  
3) The assumption A3 implies that the L\'{e}vy process $ L $ has densities 
with respect to the Lebesgue-measure on $ \R^d $ 
(see \cite{[Kol2]} p.150).\\[1ex]
The characteristic function of $ L_t $ has the following form 
(see \cite{[Sat]} p.37)
$$ \varphi_{L_t}(\xi):=\E\left[e^{i\langle\xi,L_t\rangle}\right]
          =e^{t\psi(\xi)} $$
with
$$ \psi(\xi)=\int_{\R^d\backslash\{0\}}\left(e^{i\langle\xi,x\rangle}-1-i\langle\xi,x\rangle1
\hspace{-1mm}{\rm I}_{B_1(0)}(x)\right)\nu(dx) .$$ 
We denote by $ N^L(\omega,dy,dt) $ the
Poisson-random-measure associated to the jump-process $ \Delta
L_t:=L_t-L_{t-} $ of $ L $.
The L\'{e}vy-It\^{o} decomposition theorem then gives the following path-wise
description of $ L $ 
(see \cite{[IkeWat]} p.65)
\begin{eqnarray*}
 L_t(\omega)&=& \int_0^t\int_{B_1(0)^c}yN^L(\omega,dy,ds) 
  +\int_0^t\int_{B_1(0)}y\tilde{N}^L(\omega,dy,ds),
\end{eqnarray*}
where $  \tilde{N}^L(\omega,dy,ds):=\left(N^L(\omega,dy,ds)-\nu(dy)ds\right) $
is the compensated random measure. We call $ \hat{N}^L(dy,ds):=\nu(dy)ds $ the compensator of the random measure 
$ N^L(.,dy,ds) $. Note that we do not have a drift term 
$ ct $ in the L\'{e}vy decomposition theorem since the spectral measure 
$ \mu $ is symmetric.\\[2ex]
Next we propose the assumptions on the coefficients of the SDE:\\[2mm]
{\bf A4} Let $ \Lambda\subset\R^d $ be a lattice in $ \R^d $ such that \
$ \R^d/\Lambda $ is compact.\\[1mm]
{\bf A5} Let $ \sigma:\R^d\rightarrow{\rm GL}(\R^d) $ and $
b:\R^d\rightarrow\R^d $ be $ \Lambda $-periodic and three-times continuously 
differentiable.\\[1ex]
{\bf Remarks:} Some comments on the assumptions may be helpful:\\[0.5mm]
1) The compactness in A4 is necessary to ensure that the jump-diffusion defined below
has an invariant measure.\\[0.5mm]
2) In the following we will need the existence and differentiability of the fundamental solutions for the 
jump-diffusion. Under the regularity assumptions from A5 it was proved in \cite{[Kol]} that the 
fundamental solutions are twice continuously differentiable.\\[1ex]
We consider the solution $ X $ of the following stochastic
differential equation
\begin{eqnarray*}
 &&dX_t=b(X_t)dt+\sigma(X_{t-})dL_t,\\
 && X_0=x_0.
\end{eqnarray*}
In other words, we consider the solutions of the following stochastic integral 
equation
\begin{eqnarray*}
X_t&=& x_0+\int_0^tb(X_s)ds+
\int_0^t\int_{B_1(0)^c}\sigma(X_{s-})yN^L(.,dy,ds)\\
&&+\int_0^t\int_{B_1(0)}\sigma(X_{s-})y\tilde{N}^L(.,dy,ds).\\
\end{eqnarray*}
The regularity assumtion in A5 implies the Lipschitz-continuity of the two maps 
$ x\mapsto b(x) $ and $ x\mapsto\sigma(x)y $. Existence and uniqueness of solutions to 
the above equations in the class of c\`{a}dl\`{a}g $ {\cal F}_t $-semi-martingales are well 
known facts (see \cite{[App]} p.311). 
The process $ X $ is known to be a Markov-process (see \cite{[App]} p.323). 
We denote by $ (S_t)_{t\geq0} $ the associated strongly continuous semigroup on
the space of bounded continuous functions $ C_b(\R^d) $ 
(see \cite{[App]} p.339). 
The generator of $ S $ restricted to $ C^2(\R^d)\cap C_b(\R^d) $ is the 
following integral-differential operator:
\begin{eqnarray*}
Au(x) \!&=& \!\int_{\R^d\backslash\{0\}}\left(u(x)-u(x+\sigma(x)y)
  -1\hspace{-1mm}{\rm I}_{B_1(0)}(y)
  \langle\sigma(x)y,\nabla u(x)\rangle\right)\nu(dy)\\
  \!&&\!+\langle b(x),\nabla u(x)\rangle ,
\end{eqnarray*}
(see \cite{[App]} p.340). When restricted to $ C_c^\infty(\R^d) $ the 
generator $ A $ can be regarded as a pseudo-differential operator
\begin{eqnarray*}
Au(x)=(2\pi)^{-d/2}\int_{\R^d}e^{i\langle
x,\xi\rangle}q(x,\xi)\hat{u}(\xi)d\xi
\end{eqnarray*} 
with continuous, negative definite L\'{e}vy-Khinchin-type symbol
\begin{eqnarray*}
q(x,\xi)=i\langle
  b(x),\xi\rangle+\int_{\R^d\backslash\{0\}}\left(e^{-i\langle
y,\xi\rangle}-1-i\langle y,\xi\rangle1\hspace{-1mm}{\rm
I}_{B_1(0)}(y)\right)\nu(x,dy),
\end{eqnarray*}
where $$ \nu(x,A):=\int_{\R^d\backslash\{0\}}1\hspace{-1mm}{\rm
      I}_A(\sigma(x)y)\nu(dy) $$ (see \cite{[Jac]} p.335).\\[2ex]
We will denote by $ X^\Lambda $ the image of the process $ X $ 
with respect to the canonical projection 
$ pr_\Lambda:\R^d\rightarrow\R^d/\Lambda $.
It follows from the $ \Lambda $-periodicity of the coefficients $ \sigma $ and 
$ b $, that $ X^\Lambda $ is a Markov-process on $ \R^d/\Lambda $.
The associated semigroup on $ C(\R^d/\Lambda) $ resp. generator on 
$ C^2(\R^d/\Lambda) $ will be denoted by $ (S^\Lambda_t)_{t\geq0} $ resp. 
$ A^\Lambda $. The following proposition is a version of Doeblin's celebrated 
result (see \cite{[Doo]} p.256).

\begin{proposition}\label{Spektralluecke}
There exists an invariant probability measure $ \pi $ on $ \R^d/\Lambda $ and
constants $ K,\gamma\in\R^+ $ such that for all $ f\in
C(\R^d/\Lambda) $ with $$
\overline{f}:=\int_{\R^d/\Lambda}fd\pi=0 $$ and for all $ t\geq 0 $
one has
$$ \|S_t^\Lambda f\|_{\sup}\leq Ke^{-\gamma t}\|f\|_{\sup} .$$
\end{proposition}
{\bf Proof:}
In the book of Bensoussan, Lions and Pappanicolaou one can find a version of Doeblin's result which states that for a
Markov-process with transition-probability densities on a compact manifold there exists an invariant measure,
if the transition densities are bounded from below by a constant. It is also proved that the associated semigroups 
converge exponentially as described in the present proposition (see \cite{[BenLioPap]} p.365).\\
We thus have to prove that transition-probability densities for $ X^\Lambda $ exist and are bounded from below by a 
constant $ \delta>0 $.\\[1mm]
Under the assumptions A1 to A5 Kolokoltsov proved the existence of transition-probability densities and Aronson-type
bounds for the jump-diffusion $ X $ (see \cite{[Kol]} p.749 and p.750). 
Therefore the following two facts hold:\\[1mm]
i) There exist densities $ p(t,x,y) $ with respect to the 
Lebesgue-measure for the transition-probabilities of $ X $;\\[0.5mm]
ii) There exist a $ t>0 $ such that for all compact 
$ K\subset\R^d\times\R^d $ there exists a $ \delta_0>0 $ with
$ \inf_{x,y\in K}p(t,x,y)>\delta_0 $.\\[1mm]
The densities for the transition probabilities of
$ X^\Lambda $ are then given by
$$ p_\Lambda(t,x,y)=\sum_{l\in\Lambda}p(t,x_0,l+y_0) ,$$
where $ x_0 $ resp. $ y_0 $ are arbitrary points in $ pr_\Lambda^{-1}(\{x\}) $ 
resp. $ pr_\Lambda^{-1}(\{y\}) $. 
Since $ (x,y)\mapsto p(t,x,y) $ is bounded from below by a positive constant on
each bounded subset of $ \R^d\times\R^d $, there exists a $ \delta>0 $ such 
that $ p_\Lambda(t,x,y)>\delta $ for all $ x,y\in\R^d/\Lambda $. \hfill $ \bullet $

\begin{proposition}\label{Resolvente}
For all $ f\in C(\R^d/\Lambda) $ satisfying $
\overline{f}=0 $ there exists a function $ \psi\in{\rm
Dom}(A^\Lambda) $ with $ \overline{\psi}=0 $ such that
$$ A^\Lambda\psi=f .$$
\end{proposition}
{\bf Proof:}  
By Proposition \ref{Spektralluecke} one has for all $ f $ in
$ C(\R^d/\Lambda) $ with $ \overline{f}=0 $ that  
\begin{eqnarray*}
 \left\| \int_0^\infty S_s^\Lambda fds\right\|_{\sup}
 \leq \int_0^\infty\left\|S_s^\Lambda f\right\|_{\sup}ds
 \leq \int_0^\infty K\|f\|_{\sup}e^{-\gamma s}ds<\infty.
\end{eqnarray*}
Therefore, the resolvent $$
R_\beta^\Lambda f:=\int_0^\infty e^{-\beta s}S_s^\Lambda fds $$ 
of $ A^\Lambda $ is defined in $ \beta=0 $ for all $ f\in C(\R^d/\Lambda) $ 
with $ \overline{f}=0 $ (see \cite{[Jac]} p.259).
Now, for $ \psi:=R_0^\Lambda f $ one has $ \psi\in{\rm Dom}(A^\Lambda) $ and 
$ A^\Lambda\psi=f $. From the fact that $ \overline{S_tf}=0 $ for all 
$ t\geq0 $ and the definition of $ \psi $ one sees that $ \overline{\psi}=0 $.
\hfill $ \bullet $\\[1ex]
We now come back to the periodic situation. Let $ E $ be an arbitrary set.
For all $ \Lambda $-periodic functions $ f:\R^d\rightarrow E $ there exists a 
unique $ f_\Lambda:\R^d/\Lambda\rightarrow E $ such that 
$ f=f_\Lambda\circ pr_\Lambda $.
For a periodic $ f\in C(\R^d) $ we define the mean with respect to the 
measure $ \pi $ on $ \R^d/\Lambda $ as follows
$$ \Pi(f):=\int_{\R^d/\Lambda}f_\Lambda(x)\pi(dx)=0 .$$

\begin{proposition} \label{C2}
For all $ \Lambda $-periodic $ f\in C^2(\R^d) $ satisfying $ \Pi(f)=0 $
there exists a $ \Lambda $-periodic $ \psi $ in $ C^2(\R^d) $
satisfying $ \Pi(\psi)=0 $ and
$$ A\psi=f .$$
\end{proposition}
{\bf Proof:}
Obviously, the results from Proposition \ref{Spektralluecke} and Proposition 
\ref{Resolvente} carry over to the 
restriction of $ S $ and $ A $ to $ \Lambda $-periodic functions. 
For a $ \Lambda $-periodic function $ f\in C_b(\R^d) $ with $ \Pi(f)=0 $
we then have the existence of a $ \Lambda $-periodic 
$ \psi\in{\rm Dom}(A) $ with $ \Pi(\psi)=0 $ satisfying the 
integral-differential equation $  A\psi=f $. 
The function $ \psi $ has the following representation 
$$\psi(x)=\int_0^\infty S_sf(x)ds=\int_0^\infty \int_{\R^d}p(s,x,y)f(y)dyds.$$
We still have to prove, that $ \psi $ is twice continuously differentiable.\\[1mm]
In the following we use $ \nu\circ F_x^{-1} $ to denote the image measure of
$ \nu $ with respect to the map $ F_x:y\mapsto\sigma(x)y $. 
Further, we denote by $ q_x(t,z) $ the distribution of the $ \alpha $-stable L\'{e}vy-process
with L\'{e}vy-measure $ \nu\circ F^{-1}_x $ and drift-vector $ b(x) $ at time $ t $. 
Furthermore, we will denote by $ (\partial_{z_i}q_x)(t,x-y) $ the partial differential of \
$ z\mapsto q_x(t,z) $ with respect to $ z_i $ at  $ z=x-y $.\\[1mm]   
For every multi-index $\theta $ of order smaller or equal to the degree of differentiability $ l $ of the coefficients,
V. Kolokoltsov proved that the transition densities $ p(t,x,y) $ of the jump-diffusion $ X $ satisfy the following
asymptotic representation (see \cite{[Kol]} p.743)
\begin{eqnarray*}
  (\partial_{x_i}^\theta p)(t,x,y) &=&  (\partial_{z_i}^\theta q_x)(t,x-y)\\
  && +O(t^{-l/\alpha})q_x(t,x-y)\left(\min(1,|x-y|)+t^{1/\alpha}\right).
\end{eqnarray*}
It turns out that the $ O(t^{-l/\alpha}) $-factor in the previous 
development depends only on bounds for the coefficients and its derivatives. From this follows the existence of
a bounded function $ (x,y)\mapsto C_{x,y} $ such that
\begin{eqnarray*}
  (\partial_{x_j}\partial_{x_i}p)(1,x,y)
    =(\partial_{z_j}\partial_{z_i}q_x)(1,x-y)+C_{x,y}q_x(1,x-y).
\end{eqnarray*}
Further, follows from the method developed in (see \cite{[Kol]} p.739) that 
for an appropriate constant $ C>0 $ we have
\begin{eqnarray*}
  \left|(\partial_{z_i}\partial_{z_j}q_x)(1,x-y)\right|\leq Cq_x(1,x-y). 
\end{eqnarray*}
It follows from Chapman-Kolmogorov's formula and Proposition 1 that for $ s>1 $ one has 
\begin{eqnarray*}
 \left|(\partial_{x_i}\partial_{x_j}S_sf)(x)\right|\!\!\!&\leq&
  \!\!\!\!\int_{\R^d}\!\left|(\partial_{z_i}\partial_{z_j}q_x)(1,x\!-\!y)
                \!+\!C_{x,y}q_x(1,x\!-\!y)\right|\!dy \ \|S_{s-1}f\|_{\sup}\\
  \!\!\!&\leq&\! \left(C+\sup_{x,y}|C_{x,y}|\right)\|S_{s-1}f\|_{\sup}\\
 \!\!\!&\leq& \!\left(C+\sup_{x,y}|C_{x,y}|\right)\|f\|_{\sup}e^{-\gamma(s-1)}.
\end{eqnarray*}
It is proved in the paper of Kolokoltsov that for $ t\in]0,1] $ the family of operators $ S_t $ is 
uniformly bounded in $ C^2(\R^d) $ (see \cite{[Kol]} p.749). Therefore we have that  
$ \sup_{0\leq s\leq1}\|\partial_{x_i}\partial_{x_j}S_sf\|_{\sup}<\infty $ 
if $ f\in C^2(\R^d) $.
Therefore, there exists a constant $ \tilde{C}>0 $ such that for all 
$ x\in\R^d $ one has
$$  |(\partial_{x_i}\partial_{x_j}S_sf)(x)|\leq \tilde{C}e^{-\gamma s}
 \in L^1(\R^+,ds) .$$ 
This implies that $ \psi\in C^2(\R^d) $ and
$$ (\partial_{x_i}\partial_{x_j}\psi)(x)=
   \int_0^\infty(\partial_{x_i}\partial_{x_j}S_sf)(x)ds  .$$
\hfill $ \bullet $\\[1ex]
{\bf Remark:} Let us mention here that one can find regularity results for diffusions driven 
by stable L\'{e}vy motion in \cite{[PicSav]}. However, in order to apply those 
results one has to impose the invertibility of the functions 
$ x\mapsto\sigma(x)y $ for all $ y\in\R^d $.

\section{The non-central limit theorem}

In this section we want to prove the following non-central limit theorem for 
the renormalized sequence 
$$ X^{(n)}:=n^{-1/\alpha}\left(X_{nt}-nt\Pi(b)-x_0\right) .$$ 
\begin{theorem}
Under the assumptions A1 to A5 the sequence 
$ X^{(n)} $ converges in distribution with 
respect to the Skorohod topology to an $ \alpha $-stable
L\'{e}vy-process $ X^* $ with compensator
$$ \bar{\nu}(A)=\left(\int_{\R^d/\Lambda}
 \nu\circ F_x^{-1}(A)\pi(dx)\right) ,$$
where for $ x\in\R^d/\Lambda $ we denote by $ \nu\circ F_x^{-1} $ the image 
measure of $ \nu $ with respect to the linear map 
$$ F_x:\R^d\rightarrow\R^d;y\mapsto\sigma_\Lambda(x)y. $$ 
\end{theorem}
{\bf Proof:} By Proposition \ref{C2} the equation 
$ A\psi=b-\Pi(b) $ has a $ \Lambda $-periodic solution in 
$ C^2(\R^d) $. Now, by It\^{o}'s formula one has (see \cite{[IkeWat]} p.66)
\begin{eqnarray*}
\psi(X_t)&=&\psi(X_0)+\int_0^tA\psi(X_{s})ds\\
&&+\int_0^t\int_{\R^d\backslash\{0\}}\left(\psi(X_{s-}+\sigma(X_{s-})y)
  -\psi(X_{s-})\right)\tilde{N}^L(.,dy,ds)\\
&=&\psi(X_0)+ \int_0^t(b(X_{s})-\Pi(b))ds \\
&&+\int_0^t\int_{\R^d\backslash\{0\}}\left(\psi(X_{s-}+\sigma(X_{s-})y)
  -\psi(X_{s-})\right)\tilde{N}^L(.,dy,ds).
\end{eqnarray*}
Therefore, it follows from the equation for $ X $ that
\begin{eqnarray*}
X_t-t\Pi(b)
\!\!&=&\!\! X_0+\!\int_0^t(b(X_{s})-\Pi(b))ds
+\!\int_0^t\!\int_{\R^d\backslash\{0\}}\!\!\sigma(X_{s-})y\tilde{N}^L(.,dy,ds)\\
\!\!&=&\!\!X_0+\psi(X_t)-\psi(X_0)
   +\int_0^t\int_{\R^d\backslash\{0\}}\sigma(X_{s-})y\tilde{N}^L(.,dy,ds)\\
&&-\int_0^t\int_{\R^d\backslash\{0\}}
  \left(\psi(X_{s-}+\sigma(X_{s-})y)-\psi(X_{s-})\right)
                  \tilde{N}^L(.,dy,ds)\\
\end{eqnarray*}
We can therefore decompose the sequence $ X^{(n)} $ as follows
$$  X^{(n)}_t=n^{-1/\alpha}(\psi(X_{nt})-\psi(X_0))+K^{(n)}_t+Q^{(n)}_t $$
where
$$   Q^{(n)}_t:=n^{-1/\alpha}\int_0^{nt}\int_{\R^d\backslash\{0\}}\sigma(X_{s-})y\tilde{N}^L(.,dy,ds) $$
and 
$$   K^{(n)}_t:=n^{-1/\alpha}\int_0^{nt}\int_{\R^d\backslash\{0\}}
  \left(\psi(X_{s-}+\sigma(X_{s-})y)-\psi(X_{s-})\right)
                  \tilde{N}^L(.,dy,ds) $$
are martingales. Since $ \psi $ is bounded, the term $ n^{-1/\alpha}(\psi(X_{nt})-\psi(X_0)) $ obviously 
converges toward zero as $ n\rightarrow\infty $.\\[1mm]
We will now prove that the martingales $ K^{(n)} $ converge toward zero in 
probability with respect to the Skorohod-metric. For this it is sufficient to prove the weak convergence
toward the zero-process. This follows from a general theorem which can be found in the book of Jacod and Shiryaev (see \cite{[JacShi]}). In order to apply those results we need to introduce the notions of characteristics for general semi-martingales as presented by Jacod and Shiryaev. We start with a fixed
truncation-function, which is used to seperate large jumps from smaller ones:
\[ h:\R^d\rightarrow\R^d,x\mapsto\left\{\begin{array}{cc}
     x & {\rm for} \ |x|<1 \\
     \frac{1}{|x|}x & {\rm for} \ |x|\geq1.\end{array}\right. \]
Let $ Y $ be a general semi-martingale with-respect to a filtration $ {\cal F}_t $ satisfying the usual 
conditions. The process $ \check{Y}^{(h)}_t:=\sum_{s\leq t}(\Delta Y_s-h(\Delta Y_s)) $ 
performs only the jumps of size larger than one of the process $ Y $.
The truncated process $ Y^{(h)}:=Y-\check{Y}^{(h)} $ has no jumps of size larger that one and 
can be decomposed in the following way:
$$  Y^{(h)}_t=Y_0^{(h)}+M^{(h)}_t+B^{(h)}_t ,$$
where $ B^{(h)} $ is a $ {\cal F}_t $-predictable process with finite 
variation and $ M^{(h)} $ is a local $ {\cal F}_t $-martingale 
(see \cite{[JacShi]} p.76).\\[1ex]
The predictable process $ B^{(h)} $ is called the first characteristic of 
$ Y $.\\[1mm]
For the vector-valued process $ M^{(h)} $ we denote by $ (M^{(h)})^T $ its transpose. The process
$ M^{(h)}\cdot(M^{(h)})^T $ is then a matrix-valued process.
There exists a unique $ {\cal F}_t $-previsible matrix-valued process $ \tilde{C}^{(h)} $ such that
the process $ M^{(h)}\cdot(M^{(h)})^T-\tilde{C}^{(h)} $ becomes an $ {\cal F}_t $-martingale. The process  
$ \tilde{C}^{(h)} $ is called the modified second characteristic of $ Y $.\\[1mm]
The third characteristic of $ Y $ is the $ {\cal F}_t $-predictable 
compensator $ \hat{N}^Y(.,dy,dt) $ of the random measure 
$ N^Y(.,dy,dt) $ associated to the jumps $ \Delta Y $ of $ Y $.\\[1ex]
Although the second modified characteristic and the first characteristic heavily depend on the truncation 
function $ h $ we will in the following use the symbols $ B $ and $ \tilde{C} $ to denote the first and 
modified second characteristics of $ Y $. Moreover, for a random measure $ \zeta $ on 
$ \R^d\backslash\{0\}\times[0,\infty[ $ and
a measurable function $ g:\R^d\backslash\{0\}\rightarrow\R $ we use the convenient notation
$$ (g\ast\zeta)_t:=\int_0^t\int_{\R^d\backslash\{0\}}g(y)\zeta(.,dy,ds) ,$$
whenever the random integrals are defined almost surely (see \cite{[JacShi]} 
p.66).\\[1ex]
Now let $ {\cal F}^{(n)}_t $ be a sequence of $ \sigma $-algebras and $ Y^{(n)} $ a sequence  of 
$ {\cal F}^{(n)}_t $-semi-martingales.
We will denote the three characterisics of $ Y^{(n)} $ by $ B^{(n)} $, 
$ \tilde{C}^{(n)} $ and $ \hat{N}^{(n)}(\omega,dy,dt) $.\\[1mm]
In order to prove the convergence in distribution
of $ Y^{(n)} $ toward the $ {\cal F}_t $ semimartingale $ Y $ with characteristics $ B $, 
$ \tilde{C} $ and $ \hat{N}(.,dy,dt) $ with respect to the Skorohod topology we need 
to prove (see \cite{[JacShi]} p.459):\\[1mm]  
1) $ B^{(n)} $ converges toward $ B $ with respect to sup-norm in probability;\\[0.5mm]
2) $ \tilde{C}^{(n)}_t\rightarrow \tilde{C}_t $ in probability for all 
$ t\geq0 $;\\[0.5mm]
3) $ (g\ast\hat{N}^{(n)})_t\rightarrow(g\ast\hat{N})_t $ in probability for all $ t\geq0 $ and all 
$ g\in C_b(\R^d) $ vanishing in a neighborhood of zero.\\[1ex]
We now apply this result to show that $ K^{(n)} $ converges in distribution to zero and that the
sequence $ Q^{(n)} $ converges toward the L\'{e}vy-process $ X^* $.\\[1mm]
We define the $ \sigma $-algebras $ {\cal F}^{(n)}_t:={\cal F}_{nt} $ and obtain a sequence of filtrations satisfying the usual conditions. 
The truncation of $ K^{(n)} $ with respect to $ h $ is given by
$$ K^{h,(n)}=n^{-1/\alpha}\int_0^{nt}\int_{\R^d\backslash\{0\}}
        h(\psi(X_{s-}+\sigma(X_{s-})y))-\psi(X_{s-}))\tilde{N}^L(.,dy,ds), $$
which is a martingale. Thus the first characteristic $ B^{K,(n)} $ of $ K^{(n)} $ with respect to $ h $
is always zero. 
The second modified characteristic of $ K^{(n)} $ is the predictable matrix-valued process 
$ \tilde{C}^{K,(n)} $ such that the matrix-valued process
$ K^{(n)}\cdot(K^{(n)})^T-\tilde{C}^{K,(n)} $ becomes a martingale. 
Since the local martingale $ K^{(n)} $ is a quadratic pure jump-process it 
has no continuous part in its quadratic variation-process $ C^{K,(n)} $ (see \cite{[Pro]} p.70).\\
Since there is no $ t\geq0 $ where $ \Delta L_t>0 $ with positive 
probability, there is also no $ t>0 $ such that $ \Delta K^{(n)}_t>0 $ with 
positive probability, i.e. the process $ K^{(n)} $ has no fixed time of   
discontinuity (see \cite{[JacShi]} p.101).\\[1ex] 
The second modified characteristic of $ K^{(n)} $ with respect to $ h $ is then given by 
(see \cite{[JacShi]} p.79)
\begin{eqnarray*}
 \tilde{C}^{K,(n)}_t&:=&\int_0^{nt}\int_{\R^d\backslash\{0\}}
                  H_n^T(X_{s-},y)H_n(X_{s-},y)\hat{N}^L(.,dy,ds)\\
    &=&\int_0^{nt}\int_{\R^d\backslash\{0\}}H_n^T(X_{s-},y)H_n(X_{s-},y)\nu(dy)ds
\end{eqnarray*}
where $ H_n^T(x,y) $ denotes the transpose of the vector-valued function 
$$ H_n(x,y):=h(n^{-1/\alpha}(\psi(x+\sigma(x)y))-\psi(x))) .$$ 
From this follows that
\begin{eqnarray*}
 \left|\tilde{C}^{K,(n)}_t\right|&\leq& 
   n^{1-2/\alpha}\int_0^t\int_{\R^d\backslash\{0\}}\left|G(X_{ns-},y)\right|\nu(dy)ds,
\end{eqnarray*}
with 
$$ G(x,y):=(\psi(x)+\sigma(x)y)-\psi(x))(\psi(x+\sigma(x)y)-\psi(x))^T .$$
Since $ \alpha<2 $ it follows that $ \tilde{C}^{K,(n)}_t $ converges toward zero as 
$ n\rightarrow\infty $. 
The integral of a bounded continuous function 
$ g:\R^d\rightarrow\R $ with respect to the third characteristic $ \hat{N}^{K,(n)}(\omega,dy,ds) $ of 
$ K^{(n)} $ is given by
$$  (g\ast\hat{N}^{K,(n)})_t\!=\!\!\int_0^{nt}\!\!\!\int_{\R^d\backslash\{0\}}\!\!\!
g(n^{-1/\alpha}(\psi(X_{s-}+\sigma(X_{s-})y)-\psi(X_{s-})))\hat{N}^L(.,dy,ds) .$$
Since $ \psi $ is bounded it is obvious that this converges toward zero as $ n\rightarrow\infty $. It therefore follows that all three characteristics of $ K^{(n)} $ converge toward zero. This implies the convergence of $ K^{(n)} $ toward zero in distribution. From this follows the convergence of 
$ K^{(n)} $ toward zero in probability with respect to Skorohod topology.\\[1ex]
We have proved so far that the difference $ X^{(n)}-Q^{(n)} $ converges to zero in probability with
respect to Skorohod metric. To prove the theorem it is sufficient to prove that $ Q^{(n)} $ 
converges in distribution toward the $ \alpha $-stable process $ X^* $ (see \cite{[EthKur]} p.110).
We will denote the characteristics of $ X^* $ by $ B^* $, $ \tilde{C}^* $ and 
$ \hat{N}^*(.,dy,ds) $.\\[1ex]
We now analyze the sequence of martingales $ Q^{(n)} $.
The truncation of $ Q^{(n)} $ with respect to $ h $ has the following expression
\begin{eqnarray*}
Q^{h,(n)}_t=\int_0^{nt}\int_{\R^d\backslash\{0\}}h(n^{-1/\alpha}\sigma(X_{s-})y)\tilde{N}^L(.,dy,ds).
\end{eqnarray*} 
Therefore, the process $ Q^{h,(n)} $ is a martingale and its first characteristic $ B^{Q,(n)} $ must 
vanish. Since the L\'{e}vy-process $ X^* $ also has no drift after truncation with the symmetric function 
$ h $ the convergence of the first characteristics follows.\\[1ex] 
We first compute the third characteristic and the second modified characteristic of $ Q^{(n)} $.
The third characteristic is the predictable compensator of the random measure 
$ N^{Q,(n)}(\omega,dy,dt) $ associated to the jumps $ \Delta Q^{(n)} $ of $ Q^{(n)} $. 
This is the random measure on $ \R^d\backslash\{0\}\times\R^+_0 $ which is given by 
\begin{eqnarray*}
  N^{Q,(n)}(.,A,[0,t])=\int_0^{nt}\int_{\R^d\backslash\{0\}}
 1\hspace{-1mm}{\rm I}_A(n^{-1/\alpha}\sigma(X_{s-})y))
  N^L(.,dy,ds).
\end{eqnarray*}
Its compensator is given by
\begin{eqnarray*}
 \hat{N}^{Q,(n)}(.,A,[0,t])&=&\int_0^{nt}\int_{\R^d\backslash\{0\}}
 1\hspace{-1mm}{\rm I}_A(n^{-1/\alpha}\sigma(X_{s-})y)\hat{N}^L(.,dy,ds)\\
 &=&\int_0^{nt}\int_{\R^d\backslash\{0\}}
 1\hspace{-1mm}{\rm I}_A(n^{-1/\alpha}\sigma(X_{s-})y)\nu(dy)ds.   
\end{eqnarray*}
The modified second characteristic $ \tilde{C}^{Q,(n)} $ with respect to the
truncation function $ h $ is given by the following matrix-valued process 
(see \cite{[JacShi]} p.79)
\begin{eqnarray*}
  \tilde{C}_t^{Q,(n)}
 &=&\int_0^{nt}\int_{\R^d\backslash\{0\}}(hh^T)(n^{-1/\alpha}\sigma(X_{s-})y)\hat{N}^L(.,dy,ds)\\
 &=&\int_0^{nt}\int_{\R^d\backslash\{0\}}(hh^T)(n^{-1/\alpha}\sigma(X_{s-})y)\nu(dy)ds,
\end{eqnarray*}
where $ h^T $ stands for the transpose of the vector-valued function $ h $.\\[1ex]
In order to prove the convergence of $ Q^{(n)} $ toward $ X^* $, we have to show that for all continuous
bounded functions $ g:\R^d\rightarrow\R $ vanishing in a neighbourhood of zero 
\begin{eqnarray*}
 (g\ast\hat{N}^{Q,(n)})_t=\int_0^{nt}\int_{\R^d\backslash\{0\}}g(n^{-1/\alpha}\sigma(X_{s-})y)\nu(dy)ds 
\end{eqnarray*}
converges in probability toward 
$$ (g\ast\hat{N}^*)_t=t\int_{\R^d\backslash\{0\}}g(y)\bar{\nu}(dy) $$
and that 
$ C^{Q,(n)}_t=(hh^T\ast\hat{N}^{Q,(n)})_t $ converges in probability toward 
$$  \tilde{C}^*_t=t\int_{\R^d\backslash\{0\}}hh^T(y)\bar{\nu}(dy) . $$
For this we prove that for all bounded continuous functions $ f:\R^d\rightarrow\R $ which are bounded by the 
function $ x\mapsto|x|^2 $ in a neighbourhood of zero, we have as $ n\rightarrow\infty $ that 
\begin{eqnarray*}
  (f\ast\hat{N}^{Q,(n)})_t\longrightarrow t\int_{\R^d\backslash\{0\}}f(y)\bar{\nu}(dy) \ \ {\rm in} \ L^2(\Omega,\P) .
\end{eqnarray*}
By the scaling property of $ \nu $ we can prove that
\begin{eqnarray*}
 (f\ast\hat{N}^{Q(n)})_t
 &=&n\int_0^t\int_{\R^d\backslash\{0\}}f(n^{-1/\alpha}\sigma(X_{ns-})y)\nu(dy)ds\\
 &=&    \int_0^t\int_{\R^d\backslash\{0\}}f(\sigma(X_{ns-})y)\nu(dy)ds.
\end{eqnarray*}
We will see below, that  for $ n\rightarrow\infty $ the last 
expression converges in $ L^2(\Omega,\P) $ to
\begin{eqnarray*}
  (g\ast\hat{N}^*)_t&=&  \int_0^t\int_{\R^d\backslash\{0\}}g(y)\int_{\R^d/\Lambda}
     \nu\circ F_x^{-1}(dy)\pi(dx)ds\\
&=&\int_0^t\int_{\R^d/\Lambda}\int_{\R^d\backslash\{0\}}
              g(\sigma_\Lambda(x)y)\nu(dy)\pi(dx)dt ,
\end{eqnarray*} 
where $ F_x:\R^d\rightarrow\R^d;y\mapsto \sigma_\Lambda(x)y $.\\[1mm]
The $ L^2(\Omega,\P) $-convergence follows if we can prove that for $ n\rightarrow\infty $
\begin{eqnarray*}
&&\E\left[\left(\int_0^tp(X_{ns-})ds\right)^2\right]
=2\int_0^t\int_0^s\E\left[p(X_{ns-})p(X_{nr-})\right]drds\longrightarrow0,
\end{eqnarray*}
where $$ p(z):= \int_{\R^d\backslash\{0\}}\left(f(\sigma(z)y)
-\int_{\R^d/\Lambda}f(\sigma_\Lambda(x)y)\pi(dx)\right)\nu(dy) .$$ By the
Markov-property of $ X $ and Proposition \ref{Spektralluecke} one
has
\begin{eqnarray*}
&&2\int_0^t\int_0^s\E\left[p(X_{ns-})p(X_{nr-})\right]drds\\
&=&2\int_0^t\int_0^s\E\left[\E[p(X_{ns-})|{\cal F}_{nr}]p(X_{nr-})\right]drds\\
&=&2\int_0^t\int_0^s\E\left[(S_{n(s-r)}p)(X_{nr-})p(X_{nr-})\right]drds\\
&\leq&2\int_0^t\int_0^sKe^{n(s-r)\gamma}\|p\|_{\sup}^2drds\\
&=&\frac{2K\|p\|_{\sup}^2}{n\gamma}\int_0^t(1-e^{-ns\gamma})ds
\leq\frac{2K\|p\|_{\sup}^2t}{n\gamma}\longrightarrow0.
\end{eqnarray*}
This proves that the characteristics of $ Q^{(n)} $ converge toward the characteristics of the 
L\'{e}vy-process $ X^* $. The limit-theorem from \cite{[JacShi]} described above then implies the convergence in distribution of $ Q^{(n)} $ toward $ X^* $. Since the difference of $ X^{(n)} $ and 
$ Q^{(n)} $ converges to zero in probability with respect to the Skorohod-metric the proof of the 
theorem is complete.
\hfill $ \bullet $\\[1ex]
{\bf Remark:} We note that in the case of Brownian motion, which was treated in \cite{[BenLioPap]} and
\cite{[Bha]} the martingales $ K^{(n)} $ do not converge toward zero in distribution. Since the 
scaling used in those papers is $ n^{-1/2} $ the quadratic variation does not converge toward zero and  one obtains an additional randomness in the limiting process. \\[1ex]  
{\bf Acknowledgement:} The manuscript has evolved a lot during the refereeing procedure. The author wants 
to thank all those people who contributed through their advices and questions. Also the author wants to thank Jean Picard and Vassili Kolokoltsov for their help on regularity issues.

\end{document}